\documentclass[12pt]{amsart}
\usepackage{a4,amscd,amssymb,amsfonts, epsf,graphics,caption2}
\usepackage[section]{algorithm}
\usepackage{algorithmic,color}

\usepackage[pdftex]{graphicx} % this is to include encapsulated
\usepackage{epstopdf}         % postscript files in your pdf via pdflatex

\usepackage{pstricks}

\newcommand{\supp}{\mathop{\rm supp}}

\newcommand{\R}{{\mathbb R}}

\newcommand{\C}{{\mathbb C}}

\renewcommand{\S}{\mathbb{S}}

\newcommand{\norm}[1]{\left\|#1\right\|}

\theoremstyle{plain}
\newtheorem{theorem}{Theorem}[section]

\newtheorem{lemma}[theorem]{Lemma}
\newtheorem{proposition}[theorem]{Proposition}
\theoremstyle{definition}
\newtheorem{definition}{Definition}[section]
\theoremstyle{remark}

\newtheorem*{acknowledgements}{Acknowledgements}

\numberwithin{equation}{section}

\title[Uniqueness for an inverse problem with partial data]{Uniqueness for an inverse problem in electromagnetism with partial data}
\author[B. M. Brown, M. Marletta and J. M. Reyes]{B. M. Brown$^1$, M. Marletta$^2$ and J. M. Reyes$^1$}
\date{May 12, 2015}
\address{}
\email{brownbm@cardiff.ac.uk} \email{MarlettaM@cardiff.ac.uk}
\email{reyes.juanmanuel@gmail.com}

\begin{document}

%\thanks{\emph{Keywords}: bla bla}

\maketitle

%\tableofcontents

\begin{center}
{\small $^1$School of Computer Science $\&$ Informatics,\\ Cardiff
University,\\ Cardiff CF24 3AA, United Kingdom}
\end{center}

\begin{center}
{\small $^2$School of Mathematics,\\ Cardiff University,\\
 Cardiff CF24 4AG, United Kingdom}
\end{center}

\begin{abstract}
A uniqueness result for the recovery of the electric and
magnetic coefficients in the time-harmonic Maxwell equations from
local boundary measurements is proven. No special geometrical
condition is imposed on the {\it inaccessible} part of the
boundary of the domain, apart from imposing that the boundary of
the domain is $C^{1,1}$. The coefficients are assumed to coincide
on a neighbourhood of the boundary, a natural property in
applications.

\vskip 3mm


\end{abstract}

\section{Introduction}

Let $\mu$, $\varepsilon$, $\sigma$ be positive functions on a
nonempty, bounded, open set $\Omega$ in $\R^3$, describing the
permeability, permittivity and conductivity, respectively, of an
inhomogeneous, isotropic medium $\Omega$. Let
$\partial\Omega$ denote the boundary of $\Omega$ and $N$ the
outward unit vector field normal to the boundary. Consider the
electric and magnetic fields, $E$, $H$, satisfying the so-called
time-harmonic Maxwell equations at a frequency $\omega>0$, namely
\begin{equation}\label{form:maxwelleq}
\left\{\begin{array}{ll}
\nabla\times H +i\omega \gamma E = 0, &\\
\nabla\times E -i\omega\mu H = 0, &\\
\end{array}\right.
\end{equation}
in $\Omega$, where $\gamma = \varepsilon + i\sigma /\omega$, $i$
denotes the imaginary unit, and $\nabla\times$ denotes the {\it
curl} operator.

Let $\varepsilon$, $\sigma$, $\mu$ be non-negative coefficients and assume that $\varepsilon$, $\mu$ are bounded from below in $\Omega$. Then there exist positive values of $\omega$ for which the equations
\eqref{form:maxwelleq}, posed in proper spaces and domains, with the tangential
boundary condition either $N\times H|_{\partial\Omega}=0$ or
$N\times E|_{\partial\Omega}=0$, have non-trivial solutions (see
\cite{SICh}, \cite{L}). Such values of $\omega$ are called {\it resonant
frequencies}.

The boundary data corresponding to the inverse boundary value
problem (IBVP) for the system \eqref{form:maxwelleq} only can be
given by a boundary mapping (the {\it impedance} or {\it
admittance} map) if $\omega$ is not a resonant frequency. The fact
that the position of resonant frequencies depends on the unknown
coefficients (as it is stated in \cite{OPS03}) motivated Pedro
Caro in \cite{Ca10} to consider a {\it Cauchy data set} instead of
a boundary map as boundary data. Cauchy data sets have been used
in \cite{BuU,SZ09,SZ10,Ca10,Ca11,CZ}.

This work is focused on the IBVP for the system
\eqref{form:maxwelleq} with local boundary measurements
established by a Cauchy data set taken just on a part of
$\partial\Omega$. More precisely, Definition \ref{def:domain}
describes the conditions for the domain and the part of its
boundary where the measurements are taken and Definition
\ref{def:cauchydata} (used in \cite{Ca11}) introduces the boundary
data for the IBVP studied in this article.

\begin{definition}\label{def:domain}
Let $\Omega\subset\R^3$ be a non-empty, bounded domain in $\R^3$
with $C^{1,1}$ boundary $\partial\Omega$. Assume $\Gamma$ is a
smooth proper non-empty open subset of $\partial\Omega$. We call
$\Gamma$ the {\it accessible part} of the boundary
$\partial\Omega$ and
$\Gamma_c:=\partial\Omega\setminus\overline{\Gamma}$ the {\it
inaccessible part} of the boundary.
\end{definition}

\begin{definition}\label{def:cauchydata}Let $\Omega$ and $\Gamma$
be as in Definition \ref{def:domain}. For a pair of smooth
coefficients $\mu$, $\gamma$ on $\Omega$ according to Definition
\ref{def:admissible}, define the Cauchy data set restricted to
$\Gamma$, write $C(\mu,\gamma;\Gamma)$, at frequency $\omega>0$ by
the set of couples $(T,S)\in TH_0(\Gamma)\times TH(\Gamma)$ such
that there exists a solution $(E,H)\in (H(\Omega;\mbox{curl}))^2$
 of \eqref{form:maxwelleq} in $\Omega$ satisfying $N\times E
|_{\partial\Omega} = T$ and $N\times H |_{\Gamma}=S$, where the
spaces $H(\Omega;\mbox{curl})$, $TH(\Gamma)$, $TH_0(\Gamma)$ are
defined in Definition \ref{def:spaces}.
\end{definition}

It is known that if the domain $\Omega$ is not convex and its
boundary is not $C^{1,1}$, Maxwell equations may not admit
solutions in $H^1(\Omega)$ even for boundary data in
$H^{1/2}(\partial\Omega)$ (see \cite{BS87,BS93,Sara80,Sara81,CD}.)
Thus, for a less regular domain (e.g., Lipschitz), some
non-standard Sobolev spaces are necessary. Some of them, which
will be used in these notes, appear in the following
\begin{definition}\label{def:spaces} Let $\Omega$ and $\Gamma$ be as in Definition
\ref{def:domain}. Define $ H^{1/2}(\Gamma) = \{f|_{\Gamma}\, :\, f\in
H^{1/2}(\partial\Omega)\}$, with norm $
\norm{g}_{H^{1/2}(\Gamma)}=\inf\{\norm{f}_{H^{1/2}(\partial\Omega)}\,
:\, f|_{\Gamma} =g\}$, and $ H^{1/2}_0(\Gamma) = \{f\in
H^{1/2}(\partial\Omega)\, :\, \mbox{supp}\, f\subset
\overline{\Gamma}\}$, with norm $ \norm{f}_{H^{1/2}_0(\Gamma)} =
\norm{f}_{H^{1/2}(\partial\Omega)}$. Write
$H^{-1/2}(\partial\Omega)$ for the dual space of
$H^{1/2}(\partial\Omega)$. Consider the space $H(\Omega;\mbox{curl})=\{u\in
L^2(\Omega;\C^3)\, :\, \nabla\times u\in L^2(\Omega;\C^3)\}$ with
the usual graph norm, and the following
\begin{align*}
&TH(\partial\Omega)=\{u\in H^{-1/2}(\partial\Omega;\C^3)\, :\,
N\times v|_{\partial\Omega} = u,\mbox{ for some }v\in H(\Omega;\mbox{curl})\},\\
&TH(\Gamma)=\{u |_{\Gamma}\, : \, u\in
TH(\partial\Omega)\},\\
&TH_0(\Gamma)=\{u\in TH(\partial\Omega)\, :\,
\mbox{supp}\,u\subset \overline{\Gamma}\,\},
\end{align*}
with norms
$\norm{u}_{TH(\partial\Omega)}=\inf\{\norm{v}_{H(\Omega;\mbox{curl})}\,
:\, v\in H(\Omega;\mbox{curl}),\, N\times v|_{\partial\Omega} =
u\}$, $\norm{u}_{TH(\Gamma)} =
\inf\{\norm{v}_{TH(\partial\Omega)}\, :\, v|_{\Gamma}=u\}$,
$\norm{u}_{TH_0(\Gamma)}=\norm{u}_{TH(\partial\Omega)}$. Some
properties of these spaces can be found, e.g., in \cite{Mit},
\cite{Ca10}, \cite{Ca11}.

\end{definition}

Next, the class of admissible coefficients for the uniqueness
result proven in this article is set.

\begin{definition} \label{def:admissible}Let $\Omega$ be as in Definition \ref{def:domain}. Let $M>0$. The pair of coefficients $\mu$,$\gamma$ is {\it admissible} if
$\mu,\gamma\in C^{1,1}(\overline{\Omega})\cap
W^{2,\infty}(\Omega)$, and the following conditions are satisfied:

$\bullet$ $\mbox{Re}\,\gamma\geq M^{-1}$, $\mu\geq M^{-1}$ in
$\Omega$,

$\bullet$ $\norm{\gamma}_{W^{2,\infty}(\Omega)} +
\norm{\mu}_{W^{2,\infty}(\Omega)}\leq M$.

\end{definition}

The main result of this work reads as follows.

\begin{theorem}\label{theo:main} Let $\Omega$ and $\Gamma$ be as in Definition
\ref{def:domain} and $\omega>0$ the time-harmonic frequency.
Assume $\mu_j$,$\gamma_j$ (with $j=1,2$) are two pairs of
admissible coefficients such that $\supp(\mu_1-\mu_2)$,
$\supp(\gamma_1-\gamma_2)\subset\Omega$. Then, if
$C(\mu_1,\gamma_1;\Gamma) = C(\mu_2,\gamma_2;\Gamma)$ then
$\gamma_1=\gamma_2$ and $\mu_1 = \mu_2$ in $\Omega$.
\end{theorem}

The IBVP for Maxwell equations can be seen as a vector
generalization of the inverse conductivity problem of Calder\'on.
In his seminal paper \cite{C}, Calder\'on posed two questions as
follows: Firstly, is it possible to uniquely determine the
conductivity of an unknown object from boundary measurements?
Secondly, in the affirmative case, can this conductivity be
reconstructed? Here the boundary measurements are determined by the
Dirichlet-to-Neumann map $\Lambda_{\sigma}$, which for a
conductivity $\sigma\in L^{\infty}(\Omega)$ defined on a bounded
domain $\Omega$ modelling the object, is defined by
$\Lambda_{\sigma}f=\sigma|_{\partial\Omega} (\partial u/\partial
N)|_{\partial\Omega}$, where $f\in H^{1/2}(\partial\Omega)$, $u\in
H^{1}(\Omega)$ solves the Dirichlet problem
$$
\nabla\cdot \sigma\nabla u = 0\mbox{ in }\Omega,\qquad u=f\mbox{
on }\partial\Omega,
$$
and $(\partial u/\partial N)|_{\partial\Omega}$ denotes the normal
derivative of $u$ on $\partial\Omega$.

Matti Lassas in \cite{La} proved that the boundary measurements of
the Calder\'on problem are a low-frequency limit of the boundary
data (impedance map) of the IBVP for time-harmonic Maxwell
equations under some restrictions.

Concerning the Calder\'on problem in the plane, there are three
main global uniqueness proofs giving reconstruction D-bar methods
based on complex geometrical optics (CGO) solutions: the
Schr\"odinger equation approach for twice differentiable $\sigma$
by Nachmann \cite{N95}, the first-order system approach for once
differentiable $\sigma$ by Brown and Uhlmann \cite{BU}, and the
Beltrami equation approach assuming no smoothness ($\sigma\in
L^{\infty}(\Omega)$) by Astala and P\"aiv\"arinta \cite{AsP}. The
assumption $\sigma\in L^{\infty}(\Omega)$ was the one originally
used by Calder\'on in \cite{C}. Several stability estimates have
been proven: \cite{BarBaR,BaFR,ClFR,FRo}.

In dimension $ n \geq 3 $ the best known uniqueness result for the
Calder\'on problem is due to Haberman and Tataru \cite{HT} for
continuously differentiable conductivities. A novel argument of
decay {\it in average} using Bourgain-type spaces is introduced
there. We cite some previous uniqueness results: the foundational
\cite{SU} for smooth conductivities by Sylvester and Uhlmann,
\cite{N88} where Nachman presents a reconstruction algorithm, and
\cite{B,BTo,PPaU}. Concerning conditional stability, the best
result is by the third author et al. in \cite{CGR} for
$C^{1,\varepsilon}$ conductivities on Lipschitz domains using the
method in \cite{HT}. A previous stability result was given by Heck
in \cite{He}. Roughly speaking, the method introduced by
Alessandrini in \cite{A} gives the main guidelines followed by
most stability methods for both the scalar and vector problems.

To deal with inverse problems from partial data for scalar
elliptic equations, two main approaches are found in the
literature in dimension $n>2$ (see \cite{AsLP} for $n=2$): using
Carleman estimates \cite{BU,KSjU,HeW} and using reflection
arguments \cite{I,HeW07}. This work applies to the vector case the
density argument shown in \cite{AU} for the scalar Schr\"odinger
equation by Gunther Uhlmann and Habib Ammari.

The IBVP for stationary Maxwell equations with global data was
originally proposed by Somersalo et al. in \cite{SICh}, where the
coefficients are supposed to deviate only slightly from constant
values. The same year the unique recovery of the parameters from
the scattering amplitude for $\mu$ constant was presented in
\cite{CP}, and a local uniqueness result for the IBVP from global
data was proven in \cite{SunU}. The first global determination
result for the IBVP with general coefficients $\mu$, $\gamma$ from
global boundary measurements was proven by Lassi P\"aiv\"arinta et
al. in \cite{OPS}, assuming $C^3$ smoothness on $\mu$, $\gamma$ on
$C^{1,1}$ domains. The proof is constructive. Later on, the proof
was simplified via a relation between Maxwell equations and a
matrix Helmholtz equation with a potential in \cite{OS}. Boundary
determination results appeared in \cite{Mc} and \cite{JMc} for
smooth boundaries. Chiral media were studied in \cite{Mc00}.
Stability from global boundary data was obtained in \cite{Ca10}.
Other inverse problems in electromagnetism in settings different
to the ones in this paper have been considered in
\cite{Sa,Ha,KLS,LYZ,KSU,CZ}.

%Carlos Kenig et al. proved uniqueness for the IBVP in certain
%anisotropic settings \cite{KSU}.
%
%On the other hand, Sarkola proved that the inverse scattering
%problem (ISP) can be reduced to the IBVP in \cite{Sa}. A stability
%result for the ISP supposing $\mu$ constant and as a result for
%the IBVP when the domain is a ball was given in \cite{Ha}. Such
%stability result for the IBVP was improved by P. Caro in
%\cite{Ca10}. Other inverse problems in electromagnetism have been
%considered in \cite{KLS}, \cite{LYZ}.

The uniqueness and stability issue of the IBVP for Maxwell
equations with local boundary data has been little studied. The
only works in this direction the authors are aware of are
\cite{COS}, \cite{Ca11}, where an extension of Isakov's method in
\cite{I} to Maxwell system is performed. Another extension of
methods used in the scalar case for partial data to vector systems
is in \cite{SZ10}, where uniqueness for a Dirac-type system is
proven following the ideas of \cite{KSjU}.

%%%%

Our proof uses the CGO solutions given in \cite{Ca10,Ca11} to a
matrix Schr\"odinger-type equation related to Maxwell system and a
Dirac-type system not related to Maxwell equations. As in Lemma
3.2 of \cite{KSU} we give an integral identity (Proposition
\ref{prop:intboundary}) involving a solution $Z_1$ to the
Schr\"odinger-type equation and a solution $Y_2$ to the Dirac
system. In order that such an identity holds for boundary data
restricted to $\Gamma$, the solutions have to satisfy certain
local homogeneous boundary conditions on $\Gamma_c$. In
\cite{Ca11} the solutions with such properties are constructed
from the CGO solutions following the reflection principle in
\cite{I}, arising this way the strong geometrical constraint on
$\Gamma_c$ of being plane or part of a sphere.

We manage to avoid this annoying restriction by the density
argument  given by Lemma 2 in \cite{AU} for the scalar
Schr\"odinger equation adapted to a vector Helmholtz equation
satisfied by the electric (and magnetic with different
coefficients) field related to $Z_1$ and another matrix
Schr\"odinger-type equation verified by $Y_2$. Here, unique
continuation principles for the aforementioned vector equations
and a stability estimate for the inverse of the Dirac-type
operator are required.

The density argument makes the assumption that the boundary is
$C^{1,1}$ and the coefficients coincide on a neighbourhood of the
boundary, the latter being a natural property in applications. The
rest of the proof is valid with just Lipschitz boundary.

%%%%

This article is organized as follows. In Section 2 some key matrix
equations are introduced with the novelty, compared to previous
works, that we use a Schr\"odinger-type equation satisfied by the solutions
to the Dirac-type equation $(P+W^{\ast})\widehat{Y}=0$ not related
to the Maxwell system. The CGO solutions to some of these equations
used in this article are recalled in Section 3. Section 4 is
devoted to showing an orthogonality identity involving the
potentials and solutions to matrix equations corresponding to two
couples of admissible coefficients. Two density results which are
essential in the proof are presented in Section 5. Finally, our
proof of uniqueness is expounded in Section 6 which contains a
demonstration of the bounded invertibility of a Dirac-type
operator with certain boundary conditions.

Throughout this work the following notation is used.

\vskip 1mm

\noindent{\it Notation}. Given a pair of coefficients
$(\mu_j,\gamma_j)$ for $j=1,2$, denote
$C_{\Gamma}^j:=C(\mu_j,\gamma_j;\Gamma)$, where
$C(\mu_j,\gamma_j;\Gamma)$ is defined in Definition
\ref{def:cauchydata}. For an expression like $U\cdot V$, with
$U,V$ $\C^m$-valued vector fields and $m$ a natural number,
$\cdot$ denotes the analytic extension to $\C^m$ of the Euclidean
real-inner product on $\R^m$. $I_m$ stands for the $m\times m$
identity matrix. For a matrix of complex entries $U$, the
expressions $U^t$ and $U^{\ast}$ stand for its transpose and
conjugate transpose, $\overline{U}^{\, t}$, respectively. For the
domain $\Omega$ and complex vector fields $U,V\in
H^{1}(\Omega;\C^m)$, denote
$$
(U|V)_{\Omega}:=\int_{\Omega}V^{\ast}U\,dx,\qquad
(U|V)_{\partial\Omega}:=\int_{\partial\Omega}V^{\ast}U\,ds,
$$
where $ds$ denotes the restriction of the Lebesgue measure of
$\R^3$ to $\partial\Omega$.

\begin{acknowledgements} This work is supported by the EPSRC
project EP/K024078/1. J.M.R. was also supported by the project MTM
2011-02568 Ministerio de Ciencia y Tecnolog\'ia de Espa\~na. J.M.R. wishes to thank Pedro Caro for his help over a warm meeting in ICMAT (Madrid, Spain) honouring Alberto Ruiz' 60th birthday, on Propositions 3.1 and 3.2 and uniqueness of Cauchy problems followed from unique continuation properties. The authors also thank Friedrich Gesztesy, Gerd Grubb, Hubert Kalf and William Desmond Evans for helpful discussions, as well as the referees for their careful reading of this work and useful comments.

\end{acknowledgements}

\section{The equations}

Here some differential systems related to Maxwell equations are
presented.

Let us start with the classical Schr\"odinger-type equation
approach by Petri Ola and Erkki Somersalo in \cite{OS}. Fix a
frequency $\omega>0$. Assume the coefficients $\mu$,$\gamma$ to be
in $C^{1,1}(\overline{\Omega})$ on a bounded domain with boundary
locally described by the graph of a Lipschitz function. Following
the notation in \cite{Ca10}, write$$ \alpha :=\log\gamma,\qquad
\beta:=\log\mu,\qquad \kappa:=\omega\mu^{1/2}\gamma^{1/2}.
$$

The vector fields $E,H\in H(\Omega;\mbox{curl})$ solve Maxwell equations
\eqref{form:maxwelleq} with coefficients $\mu,\gamma$ if and only
if $ X = (h\,\,\,\,\, H^t\,\,|\,\,\, e\,\,\,\,\,E^t)^t $ solves
the so-called {\it augmented system} $(P+V)X=0$ and the scalar
fields $e,h$ vanish, where
\begin{equation}\label{form:operatorP}
P:=\left(%
\begin{array}{cc|cc}
   &  &  & D\cdot \\
   &  & D & -D\times \\
   \hline
   & D\cdot &  &  \\
  D & D\times &  &  \\
\end{array}%
\right)\begin{array}{c}
\phantom{.}\\
\phantom{.}\\
\phantom{.}\\
,\\
\end{array}\, V:=\left(%
\begin{array}{cc|cc}
   \omega\mu &  &  & D\alpha\cdot \\
   & \omega\mu I_3 & D\alpha & \\
   \hline
   & D\beta\cdot & \omega\gamma &  \\
  D\beta &  &  & \omega\gamma I_3 \\
\end{array}%
\right) \begin{array}{c}
\phantom{.}\\
\phantom{.}\\
\phantom{.}\\
,\\
\end{array}
\end{equation}
with $D:=1/i\nabla$, $\nabla$ denoting the {\it gradient}
operator, and $\nabla\cdot$ the {\it divergence} operator.
Further, $X$ solves $(P+V)X=0$ if and only if $
Y=\mbox{diag}(\mu^{1/2}I_4\, ,\,\gamma^{1/2}I_4)\, X $ solves the
{\it rescaled system} $(P+W)Y=0$ with
\begin{equation}\label{form:1}
W:=\kappa I_8 + {1\over 2}\left(%
\begin{array}{cc|cc}
   &  &  & D\alpha\cdot \\
   &  & D\alpha & D\alpha\times \\
   \hline
   & D\beta\cdot &  &  \\
  D\beta & -D\beta\times &  &  \\
\end{array}%
\right)\begin{array}{c}
\phantom{.}\\
\phantom{.}\\
\phantom{.}\\
,\\
\end{array}
\end{equation}
where $\times$ denotes the cross product. Define the terms
\begin{equation}\label{form:matpoten}
Q:= WP-P\, W^t-WW^t,\qquad
\widehat{Q}:=W^{\ast}P-P\,\overline{W}-W^{\ast} \,\overline{W}.
\end{equation}It can be checked that the expressions $WP-P\, W^t$ and
$W^{\ast}P-P\,\overline{W}$ are zeroth-order. Therefore, the
second-order operators$$-\Delta I_8 + Q = (P+W)(P-W^t),\qquad
-\Delta I_8 + \widehat{Q}= (P+W^{\ast})(P-\overline{W}),$$which do
not contain first-order terms, are Schr\"odinger-type. In
\cite{Ca10},\cite{Ca11} a further zeroth-order operator $Q'$ is
considered which it is not used in this work.

Note that if $(-\Delta I_8 + Q )Z=0$ and
$X:=\mbox{diag}(\mu^{-1/2}I_4\, ,\,\gamma^{-1/2}I_4)\,(P-W^t)Z$
then $(P+V)X=0$. If additionally the scalar fields in $X$ are
identically zero, the vector fields in $X$ give the
electromagnetic fields verifying Maxwell equations. Moreover, if
$\widehat{Z}$ solves $(-\Delta I_8 + \widehat{Q})\widehat{Z}=0$
then $\widehat{Y}:= (P-\overline{W})\widehat{Z}$ is solution to
$(P+W^{\ast})\widehat{Y}=0$.

Finally, a new matrix Schr\"odinger potential is introduced,
namely $\widetilde{Q}$ in \eqref{form:matpotentilde}, satisfying
Lemma \ref{lemma:schro} below.

Specify explicitly the $(\kappa,\alpha,\beta)$-dependence of $W$,
$W^t$ and $W^{\ast}$ by writing
$$
W = W(\kappa,\alpha,\beta),\qquad W^t =
W^t(\kappa,\alpha,\beta),\qquad W^{\ast} =
W^{\ast}(\kappa,\alpha,\beta).
$$
A straightforward computation gives $
W^{\ast}(\kappa,\alpha,\beta) =
-W^t(-\overline{\kappa},\overline{\alpha},\beta )$. Since the
relation among $\kappa$, $\alpha$ and $\beta$ is not involved in
the proof of the fact that $-P W^t + W P$ is zeroth-order, we
deduce that $ -P W^t(-\overline{\kappa},\overline{\alpha},\beta) +
W(-\overline{\kappa},\overline{\alpha},\beta) P $ is zeroth-order.
Thus, the matrix operator
\begin{equation*}
(P+W(-\overline{\kappa},\overline{\alpha},\beta))(P-W^t(-\overline{\kappa},\overline{\alpha},\beta))=-\Delta
I_8 + \widetilde{Q}
\end{equation*}
is Schr\"odinger-type, where
\begin{equation}\label{form:matpotentilde}
\widetilde{Q}:=-PW^t(-\overline{\kappa},\overline{\alpha},\beta)+W(-\overline{\kappa},\overline{\alpha},\beta)P-W(-\overline{\kappa},\overline{\alpha},\beta)
W^t(-\overline{\kappa},\overline{\alpha},\beta)
\end{equation}
is zeroth-order. We deduce the following
\begin{lemma}\label{lemma:schro}
Assume $\widehat{Y}\in H^1(\Omega;\C^8)$. If $(P+W^{\ast})\widehat{Y}=0$ then $(-\Delta I_8 +
\widetilde{Q}\, )\widehat{Y}=0$.
\end{lemma}

\noindent {\it Notation}. In the rest of the manuscript, for two
pairs of coefficients $\mu_j$, $\gamma_j$ with $j=1,2$, we will
write $Q_j$, $\widetilde{Q}_j$, $W_j$ to refer to the zeroth-order
matrix operators $Q$, $\widetilde{Q}$, $W$ defined in
\eqref{form:matpoten},\eqref{form:matpotentilde},\eqref{form:1},
respectively, for the case $\mu=\mu_j$, $\gamma=\gamma_j$.

\section{The special solutions}\label{section:solutions}

In this section we recall the almost exponentially growing
solutions $Z$, $Y$ constructed in \cite{Ca10} for the systems
$(-\Delta I_8 + Q)Z=0$, $(P+W^{\ast})Y=0$ based on ideas of the
papers \cite{SU}, \cite{B}, \cite{OS}, \cite{KSU}. Here the
coefficients $\mu_j$, $\gamma_j$ ($j=1,2$) under Theorem
\ref{theo:main}'s conditions have to be considered extended to the
whole Euclidean space $\R^3$. We denote the extended coefficients
in the same manner $\mu_j$, $\gamma_j$. The extensions fulfill the
properties as follows:
\begin{enumerate}
\item[1.] They are Whitney type (see \cite{S} for their
construction). \item[2.] The extensions preserve the regularity
and a priori conditions stated in Definition \ref{def:admissible}.
\item[3.] The extended functions $\mu_j$, $\gamma_j$ satisfy that
$\mu_j=\mu_0$, $\gamma_j=\varepsilon_0$ outside a ball
$B(\mbox{O},\rho)$ centered at origin $\mbox{O}$ with radius
$\rho>0$ such that $\overline{\Omega}\subset B(\mbox{O},\rho)$,
where $\mu_0$, $\varepsilon_0>0$ are constants. \item[4.] Denoting
likewise the matrices $Q$, $\widehat{Q}$ obtained by replacing the
coefficients by their extensions, the matrix functions
$\omega^2\varepsilon_0\mu_0 I_8 + Q$, $\omega^2\varepsilon_0\mu_0
I_8 + \widehat{Q}$ are compactly supported in
$\overline{B(\mbox{O},\rho)}$.
\end{enumerate}

These properties allow Caro to prove Proposition 9 and Proposition
11 in \cite{Ca10}, which we present here for  the reader's
convenience.

\noindent\textit{Notation}: $ \norm{f}_{L^2_{\delta}}^2 =
\int_{\R^3}(1+|x|^2)^{\delta}|f(x)|^2\,dx $, for $\delta\in\R$.

\begin{proposition}\label{prop:special1}Let $-1<\delta<0$ and $\zeta\in\C^3$ with $\zeta\cdot\zeta = \omega^2\varepsilon_0
\mu_0$. Assume
\begin{align*}
&|\zeta|>C(\delta , \rho)\left(
\sum_{j=1,2}\norm{\omega^2\varepsilon_0\mu_0 +
q_{j}}_{L^{\infty}(B_{\rho})} + \norm{\omega^2\varepsilon_0\mu_0
I_8 + Q}_{L^{\infty}(B_{\rho};\mathcal{M}_{8\times 8})} \right),
\end{align*}where $\mathcal{M}_{8\times8}$ denotes the space of $8\times8$ matrices with complex entries, and
$$
q_{1} = -{1\over 2}\, \Delta\beta-\kappa^2-{1\over 4}\,
(D\beta\cdot D\beta),\qquad q_{2} = -{1\over 2}\,
\Delta\alpha-\kappa^2-{1\over 4}\, (D\alpha\cdot D\alpha).
$$

Then there exists a solution
$$Z(x,\zeta)=e^{i\zeta\cdot x}(L(\zeta)+R(x,\zeta)),$$to $(-\Delta I_8 +
Q)Z=0$ in $\R^3$ with $Z|_{\Omega}\in H^2(\Omega;\C^8)$, where
\begin{equation}\label{form:boundedterm}
L(\zeta)={1\over|\zeta|}\left(
\begin{array}{c}
\zeta\cdot A_1\\
\omega\varepsilon_0^{1/2}\mu_0^{1/2}B_1\\
\hline
\zeta\cdot B_1\\
\omega\varepsilon_0^{1/2}\mu_0^{1/2}A_1\\
\end{array}
\right)
\end{equation}
with $A_1$,$B_1$ constant complex vector fields, and
\begin{align}
\label{form:resR}&\norm{R}_{L^2_{\delta}(\R^3;\C^8) }\leq {C(\delta ,
\rho)\over |\zeta|}\,|L|\,\norm{\omega^2\epsilon_0\mu_0 I_8 +
Q}_{L^{\infty}(B_{\rho};\mathcal{M}_{8\times 8})}.
\end{align}

Furthermore, $Y := (P-W^t)Z$ solves $(P+W)Y = 0$ in $\R^3$ and has
the form
$$
Y = (0\,\,\,\,\,\mu^{1/2} H^t\,\,|
\,\,0\,\,\,\,\,\gamma^{1/2}E^t)^t,
$$
with $E,H$ solutions of \eqref{form:maxwelleq} in $\R^3$.

\end{proposition}

\begin{proposition}\label{prop:special2}
Let $\zeta\in\C^3$ with $\zeta\cdot\zeta = \omega^2 \varepsilon_0
\mu_0$ and
\begin{align*}
&|\zeta|>C(\rho) \norm{\omega^2\varepsilon_0\mu_0 I_8 +
\widehat{Q}}_{L^{\infty}(B_{\rho};\mathcal{M}_{8\times 8})}.
\end{align*}
Then there exists a solution
$$\widehat{Y}(x,\zeta) = e^{i\zeta\cdot x}(M(\zeta)+S(x,\zeta))$$to $(P+W^{\ast})\widehat{Y} =0$ in $\R^3$ with $\widehat{Y}|_{\Omega}\in H^1(\Omega;\C^8)$,
\begin{equation}\label{form:boundedterm2}
M(\zeta)={1\over|\zeta|}\left(
\begin{array}{c}
\zeta\cdot A_2\\
-\zeta\times A_2\\
\hline
\zeta\cdot B_2\\
\zeta\times B_2\\
\end{array}
\right),
\end{equation}
where $A_2$,$B_2$ are constant complex vector fields, and
\begin{align}
\label{form:resS}&\norm{S}_{L^2(\Omega;\C^8)} \leq {C(\rho ,
\Omega)\over |\zeta|}\,\bigg(\norm{\omega^2\varepsilon_0\mu_0 I_8
+ \widehat{Q}}_{L^{\infty}(B_{\rho};\mathcal{M}_{8\times 8})} +
\norm{W}_{L^{\infty}(\Omega;\mathcal{M}_{8\times 8})} \bigg).
\end{align}
\end{proposition}

\section{An orthogonality identity}

This section is aimed at proving an orthogonality identity given
by Proposition \ref{prop:intboundary} involving solutions on the
open set $\Omega$ to certain matrix partial differential equations
whose traces contain information supported on $\overline{\Gamma}$.

\begin{lemma}\label{lemma:lemma3.3}
Let $\mu_j$, $\gamma_j$ belong to $C^{0,1}(\overline{\Omega})$ for
$j=1,2$. Let $$ Y_1 = (0\,\,\,\,\,\mu_1^{1/2}
H_1^t\,|\,0\,\,\,\,\,\gamma_1^{1/2}E_1^t)^t,
$$
with $E_1,H_1\in H(\Omega;\mbox{curl})$ solutions of
\begin{gather}\label{form:maxwell1a}
\nabla\times H_1 +i\omega \gamma_1 E_1 = 0,\\
\label{form:maxwell1b}\nabla\times E_1 -i\omega\mu_1 H_1 = 0,
\end{gather}
in $\Omega$ such that $N\times E_1 = 0$ on $\Gamma_c$. In
addition, suppose that $$Y_2 = (f^1\,\,\,\,\, (u^1)^{t}\,|\,
f^2\,\,\,\,\, (u^2)^t)^t$$ is a solution to $(P+W_2^{\ast})Y_2 =
0$ in $\Omega$ with $f^j\in H^1(\Omega)$, $u^j\in
H(\Omega;\mbox{curl})$ and $ f^1 = N\times u^2=0 $ on $\Gamma_c$.
Hence, for any pair $E_2$, $H_2$ in $H(\Omega;\mbox{curl})$ of
solutions to
\begin{equation}\label{form:maxwell2} \left\{
\begin{array}{ll}
\nabla\times H_2 +i\omega \gamma_2 E_2 = 0,&\\
\nabla\times E_2 -i\omega\mu_2 H_2 = 0,&
\end{array}
\right.
\end{equation}
in $\Omega$ such that $N\times E_2|_{\partial\Omega}=0$ on
$\Gamma_c$, the following estimate holds:
\begin{gather*}
|(Y_1|PY_2)_{\Omega}-(PY_1|Y_2)_{\Omega}|\\
\leq C\, \bigg(\norm{N\times
(E_1-E_2)|_{\partial\Omega}}_{TH_0(\Gamma)} +
\norm{N\times(H_1-H_2)|_{\Gamma}}_{TH(\Gamma)}\bigg)\\
\times \left( \norm{\mu_2^{-1/2}}_{C^{0,1}(\overline{\Gamma})}
\norm{f^2}_{H^{1/2}(\Gamma)} +
\norm{\gamma_2^{1/2}}_{C^{0,1}(\overline{\Gamma})} \norm{N\times
u^1}_{TH(\Gamma)}\right.\\
\left. + \norm{\gamma_2^{-1/2}}_{C^{0,1}(\overline{\Gamma})}
\norm{f^1}_{H^{1/2}_0(\Gamma)} +
\norm{\mu_2^{1/2}}_{C^{0,1}(\overline{\Gamma})} \norm{N\times
u^2}_{TH_0(\Gamma)}\right)\\
+ C\, \bigg(\norm{N\times E_1}_{TH_0(\Gamma)} +
\norm{N\times H_1}_{TH(\Gamma)}\bigg)\\
\times \left(
\norm{\mu_1^{-1/2}-\mu_2^{-1/2}}_{C^{0,1}(\overline{\Gamma})}
\norm{f^2}_{H^{1/2}(\Gamma)} +
\norm{\gamma_1^{1/2}-\gamma_2^{1/2}}_{C^{0,1}(\overline{\Gamma})}
\norm{N\times
u^1}_{TH(\Gamma)}\right.\\
\left. +
\norm{\gamma_1^{-1/2}-\gamma_2^{-1/2}}_{C^{0,1}(\overline{\Gamma})}
\norm{f^1}_{H^{1/2}_0(\Gamma)} +
\norm{\mu_1^{1/2}-\mu_2^{1/2}}_{C^{0,1}(\overline{\Gamma})}
\norm{N\times u^2}_{TH_0(\Gamma)}\right).
\end{gather*}

\end{lemma}

\vskip 5mm

Lemma \ref{lemma:lemma3.3} follows from the proof of Lemma 3.3 in
\cite{Ca11} by making the solutions $Y_1$, $Y_2$ on $\Omega$ play
the role of $\mathcal{Y}_1$, $\mathcal{Y}_2$ on $U$ in Lemma 3.3
from \cite{Ca11}. This is achieved imposing directly to $Y_1$,
$Y_2$ the appropriate boundary conditions on $\partial\Omega$,
namely the tangential component of the electric field appearing in
the structure of $Y_1$ vanishes on the inaccessible part of the
boundary, and concerning $Y_2$, the trace of the first component
and the tangential component of the second vector field also
vanish on the inaccessible part of the boundary. In Lemma 3.3 from
\cite{Ca11} such boundary conditions for $\mathcal{Y}_1$,
$\mathcal{Y}_2$ come from a reflection argument using the special
geometric conditions assumed to $\partial U\setminus
\overline{\Gamma}$ there, which can not be used here.

The proof of Lemma \ref{lemma:lemma3.3} uses Lemma 2.2, Lemma 2.4,
Lemma 2.5 and Lemma 2.6 in \cite{Ca11}. Lemma
\ref{lemma:lemma3.3}'s proof is omitted since, up to these
comments, is identical to Lemma 3.3's proof in \cite{Ca11}.

\begin{proposition}\label{prop:intboundary}
Let $\mu_j$, $\gamma_j$ be an admissible pair of coefficients
($j=1,2$) such that $C_{\Gamma}^1=C_{\Gamma}^2$. Additionally,
suppose $\mu_1=\mu_2$, $\gamma_1=\gamma_2$,
$\partial_{x_l}\mu_1=\partial_{x_l}\mu_2$ and
$\partial_{x_l}\gamma_1=\partial_{x_l}\gamma_2$ on
$\overline{\Gamma}$ for $l=1,2,3$. Then
$$((Q_1-Q_2)Z_1|Y_2)_{\Omega}=0$$holds for any $Z_1\in H^1(\Omega;\C^8)$
solving $(P-W_1^t)Z_1=Y_1$ in $\Omega$ with
$$ Y_1 = (0\,\,\,\,\,\mu_1^{1/2}
H_1^t\,|\,0\,\,\,\,\,\gamma_1^{1/2}E_1^t)^t,
$$
where $E_1,H_1\in H(\Omega;\mbox{curl})$ are solutions of
\begin{gather}\label{form:maxwell1a}
\nabla\times H_1 +i\omega \gamma_1 E_1 = 0,\\
\label{form:maxwell1b}\nabla\times E_1 -i\omega\mu_1 H_1 = 0,
\end{gather}
in $\Omega$ and $N\times E_1|_{\partial\Omega} = 0$ on $\Gamma_c$,
and for any $Y_2\in H^1(\Omega;\C^8)$ verifying $(P+W_2^{\ast})Y_2 = 0$
in $\Omega$ such that $ Y_2|_{\partial\Omega}=0 $ on $\Gamma_c$.
\end{proposition}

\vskip 5mm

\noindent\textbf{Proof of Proposition \ref{prop:intboundary}}.
Following the proof of Proposition 1 in \cite{Ca11}, we obtain the
identity
\begin{gather*}
((Q_1-Q_2)Z_1|Y_2)_{\Omega} = ((W_1^t-W_2^t)Z_1|P_N
Y_2)_{\partial\Omega} + (Y_1|PY_2)_{\Omega}-(PY_1|Y_2)_{\Omega},
\end{gather*}
whose proof involves integration by parts and the relations
$(P+W_2^{\ast})Y_2=0$, $Y_1=(P-W_1^t)Z_1$, $(P+W_1)Y_1=0$. Here,
denote$$
P_N:={1\over i}\,\left(%
\begin{array}{cc|cc}
   &  &  & N\cdot \\
   &  & N & -N\times \\
   \hline
   & N\cdot &  &  \\
  N & N\times &  &  \\
\end{array}%
\right)\begin{array}{c}
\phantom{.}\\
\phantom{.}\\
\phantom{.}\\
.\\
\end{array}
$$
Since $\mu_1=\mu_2$, $\gamma_1=\gamma_2$,
$\partial_{x_l}\mu_1=\partial_{x_l}\mu_2$ and
$\partial_{x_l}\gamma_1=\partial_{x_l}\gamma_2$ on
$\overline{\Gamma}$ (for $l=1,2,3$), it follows that
\begin{equation}\label{form:25}
W_1^t-W_2^t=0\qquad \mbox{on }\overline{\Gamma}.
\end{equation}
From the fact that $Y_2|_{\partial\Omega}$ is supported on
$\overline{\Gamma}$, we have
\begin{equation}\label{form:26}
P_NY_2|_{\Gamma_c}\equiv 0.
\end{equation}

By \eqref{form:25},\eqref{form:26},
\begin{gather}\label{form:24}
((W_1^t-W_2^t)Z_1|P_N Y_2)_{\partial\Omega} =0.
\end{gather}
Since $\mu_1=\mu_2$ and $\gamma_1=\gamma_2$ on
$\overline{\Gamma}$,
\begin{gather}
\label{form:20}\norm{\mu_1^{1/2}-\mu_2^{1/2}}_{C^{0,1}(\overline{\Gamma})}=
\norm{\mu_1^{-1/2}-\mu_2^{-1/2}}_{C^{0,1}(\overline{\Gamma})}
=\norm{\gamma_1^{1/2}-\gamma_2^{1/2}}_{C^{0,1}(\overline{\Gamma})}\\
\label{form:21}=
\norm{\gamma_1^{-1/2}-\gamma_2^{-1/2}}_{C^{0,1}(\overline{\Gamma})}=0.
\end{gather}

Since $(N\times E_1|_{\partial\Omega},N\times H_1|_{\Gamma})\in
C_{\Gamma}^1=C_{\Gamma}^2$, there exist solutions $E_2$,$H_2$ in
$H(\Omega;\mbox{curl})$ to the Maxwell system
\eqref{form:maxwell2} such that $N\times E_2|_{\partial\Omega}$ is
supported on $\overline{\Gamma}$ and
\begin{equation}\label{form:22}
(N\times E_1|_{\partial\Omega},N\times H_1|_{\Gamma})=(N\times
E_2|_{\partial\Omega},N\times H_2|_{\Gamma}).
\end{equation}

Applying Lemma \ref{lemma:lemma3.3} to $Y_1$, $Y_2$, $E_2$, $H_2$
and by \eqref{form:20}-\eqref{form:21}, \eqref{form:22},
\begin{equation}\label{form:23}
(Y_1|PY_2)_{\Omega}-(PY_1|Y_2)_{\Omega}=0.
\end{equation}
By identities \eqref{form:24},\eqref{form:23}, we conclude
$((Q_1-Q_2)Z_1|Y_2)_{\Omega}=0$.

\begin{flushright}
$\square$\end{flushright}

\section{Density and unique continuation results}

In the remainder of the paper, let $\Omega$ and $\Gamma$ be as in
Definition \ref{def:domain}. If $E$, $H\in H(\Omega;\mbox{curl})$
solve system \eqref{form:maxwelleq} in $\Omega$ for certain
coefficients $\mu$,$\gamma$, then $E$, $H$ are also solutions to
the following second order system:
\begin{equation*}
\begin{array}{ll}
\nabla\times(\mu^{-1}\nabla\times E)-\omega^2 \gamma E = 0, &\\
\nabla\times (\gamma^{-1}\nabla\times H)-\omega^2\mu H=0. &\\
\end{array}
\end{equation*}

\noindent{\it Notation}. For known coefficients $\mu$,$\gamma$ and
frequency $\omega>0$, notation $L$ will refer to the following
Helmholtz-type vector second order differential operator defined
in the sense of distributions for $U\in (C^{\infty}_0)'(\Omega;\C^3)$ by
\begin{equation}\label{form:secondorder}
L \,U = \nabla\times(\mu^{-1}\nabla\times U)-\omega^2 \gamma U.
\end{equation}

For Lipschitz continuous functions $\mu$, $\gamma$ on $\Omega$, assuming $\mu$ to be
 bounded from below, $L \,U$ is an $L^2$ vector field if $U\in H(\Omega;\mbox{curl})$ and $\nabla\times(\nabla\times U)\in L^2(\Omega;\C^3)$, since
 $
 \nabla\times(\mu^{-1}\nabla\times U) = (\nabla \mu^{-1})\times(\nabla\times U)+\mu^{-1}\nabla\times(\nabla\times U).
 $

Note the following integration by parts formula for any $E,F\in C^{\infty}(\Omega;\C^3)$:
\begin{equation}\label{form:intpartsL}
\begin{array}{ll}
\int_{\Omega} (L\, E)\cdot F\, dx = \int_{\Omega} E\cdot (L\, F)\,
dx& \\ & \\- \int_{\partial\Omega} \mu^{-1}\,
(E\cdot(N\times(\nabla\times F)) + (\nabla\times E)\cdot (N\times
F))\, ds.&\\ \end{array}
\end{equation}

Next, the density result for the scalar Schr\"odinger equation
given by Lemma 2 in \cite{AU} is adapted to Schr\"odinger-type
matrix equations (Proposition \ref{prop:density}) and the second
order operator $L$ (Proposition \ref{prop:density_curlcurl}).

\begin{proposition}\label{prop:density}Let $\Omega'$ be an open subset of $\R^3$ with $C^2$ boundary. Assume $\Omega'\subset\subset\Omega$ and
$\Omega\setminus\overline{\Omega'}$ is connected. Let
$\widetilde{Q}$ be the zeroth-order $8\times 8$ matrix operator
defined in \eqref{form:matpotentilde} for an admissible pair of
coefficients $\mu$,$\gamma$. Then the set
$$
\widetilde{K}(\Omega):=\{g\in H^2(\Omega;\C^8)\, :\, (-\Delta I_8 +
\widetilde{Q})g=0\mbox{ in }\Omega,\,\, g|_{\partial\Omega}=0
\mbox{ on } \Gamma_c\}
$$
is dense in the space $ K(\Omega):=\{v\in H^2(\Omega;\C^8)\, :\,
(-\Delta I_8 + \widetilde{Q})v=0\mbox{ in }\Omega\} $ with respect
to the topology in $L^2(\Omega';\C^8)$. Here, $g$,$v$ denote $8\times 1$
vector fields.
\end{proposition}

\begin{proposition}\label{prop:density_curlcurl}Let $\Omega'$ be an open subset of $\R^3$ with $C^2$ boundary. Assume $\Omega'\subset\subset\Omega$ and
$\Omega\setminus\overline{\Omega'}$ is connected. Let $L$ be the
differential operator defined in \eqref{form:secondorder} for an admissible pair of
coefficients $\mu$,$\gamma$. Then
the set $\widetilde{N}(\Omega)$ of vector functions $\widetilde{E}\in H(\Omega;\mbox{curl})$ such that $
\nabla\times(\nabla\times \widetilde{E})\in L^2(\Omega;\C^3)$, $L\, \widetilde{E}=0$ in $\Omega$, $N\times
\widetilde{E}|_{\partial\Omega}=0$ on $\Gamma_c$, is dense in the space

\begin{center}
$ N(\Omega):=\{E\in H(\Omega;\mbox{curl})\, :\,
\nabla\times (\nabla\times E)\in L^2(\Omega;\C^3),\, L\, E =0\mbox{ in }\Omega\} $
\end{center}

\noindent with respect to the topology in $L^2(\Omega';\C^3)$.
\end{proposition}

The following unique continuation principles, Lemma
\ref{lemma:uniquecont} and Lemma \ref{lemma:uniquecontL}, are used
to prove Proposition \ref{prop:density} and Proposition
\ref{prop:density_curlcurl}.

\begin{lemma}\label{lemma:uniquecont}(\textbf{Unique continuation principle for matrix Schr\"odinger-type equations})
Let $\widetilde{Q}$ be the zeroth-order $8\times 8$ matrix operator defined in \eqref{form:matpotentilde} for a pair $\mu$, $\gamma$ of admissible coefficients. Assume $\Omega'\subset\subset\Omega$ and
$\Omega\setminus\overline{\Omega'}$ is connected with
$\partial\Omega'\in C^2$. Hence,

i) If $u\in H^2(\Omega;\C^8)$ satisfies $(-\Delta I_8 + \widetilde{Q})u = 0$ in $\Omega$ and $u=0$ on $B$ for some open ball $B$ such that $\overline{B}\subset \Omega$ then $u=0$ in $\Omega$.

ii) Suppose $u\in H^2(\Omega\setminus\overline{\Omega'}\, ;\C^8)$ verifies $(-\Delta I_8 + \widetilde{Q})u = 0$ in
$\Omega\setminus\overline{\Omega'}$, $ u=0\mbox{ on
}\partial\Omega$, $\big({\partial\over\partial\nu}\,
I_8\big)u\big|_{\partial\Omega}=0\mbox{ on }\Gamma, $ where
$\Gamma$ is a smooth proper non-empty open subset of
$\partial\Omega$. Then $u=0$ on $\Omega\setminus\Omega'$.

\end{lemma}

\begin{lemma}\label{lemma:uniquecontL}(\textbf{Unique continuation principle for $L$})
Let $\mathcal{G}$ be a nonempty, open, bounded, connected subset of $\R^3$ with Lipschitz boundary $\partial\mathcal{G}$. Let $L$ denote the operator \eqref{form:secondorder} for scalar functions $\mu$,$\gamma\in C^1(\,\overline{\mathcal{G}}\,)$ with $\mu\geq C$, $\mbox{Re}\,\gamma\geq C$  in $\mathcal{G}$ for some constant $C>0$.

Further, assume $U\in H(\mathcal{G};\mbox{curl})$ and $\nabla\times(\nabla\times U)\in L^2(\mathcal{G};\C^3)$. Therefore:

i) If $L\, U =0$ in $\mathcal{G}$ and $U=0$ in $B$ for some open ball $B$ such that $\overline{B}\subset\mathcal{G}$, then $U=0$ in $\mathcal{G}$.

ii) Let $\Gamma'$ denote a nonempty, smooth, open subset of $\partial\mathcal{G}$. If $L\, U =0$ in $\mathcal{G}$, $N\times U|_{\partial\mathcal{G}}=0$ on $\Gamma'$ and $N\times(\nabla\times U)|_{\partial\mathcal{G}}=0$ on $\Gamma'$, then $U=0$ in $\mathcal{G}$. Here $N$ also denotes the outward unit vector field normal to $\partial\mathcal{G}$.

\end{lemma}

\begin{figure}[!htp]
\centering
\begin{pspicture}(8,4.4)(4,-2.6)
%\psdiamond(1.5,1.5)(1.5,0.5)
%\pstriangle[fillstyle=solid,fillcolor=yellow](2,0)(3,1)
\pscurve(2,0)(2.1,1.5)(4,2)
\pscurve(8,2)(9.5,1.5)(9.5,0)(11,-0.8)(4,-2.1)(2,-2.1)(2,0)
\pscurve[linewidth=0.08, arrows=<->](4,2)(6,2.3)(8,2)
\pscurve(4.8,2.2)(5,3.5)(7,3.5)(7.3,2.2)
\psccurve(4,-1)(6,0)(8,-1)(6.5,-1.5) \rput(6,1.9){\small $\Gamma$}
\rput(6,0.8){\LARGE $\Omega\setminus\overline{\Omega'}$}
\rput(6,-0.8){\LARGE $\Omega'$} \rput(6,3){\Large
$\Gamma_{\varepsilon}$}
\end{pspicture}
\caption{\label{fig}Picture of possible choices of the sets
$\Gamma$, $\Gamma_{\varepsilon}$, $\Omega$, $\Omega'$ in the
plane. This is for the sake of clarification only; remember that
the open sets $\Gamma_{\varepsilon}$, $\Omega$, $\Omega'$ are
taken in the Euclidean topology of $\R^3$ and $\Gamma\subset
\partial\Omega$.}
\end{figure}
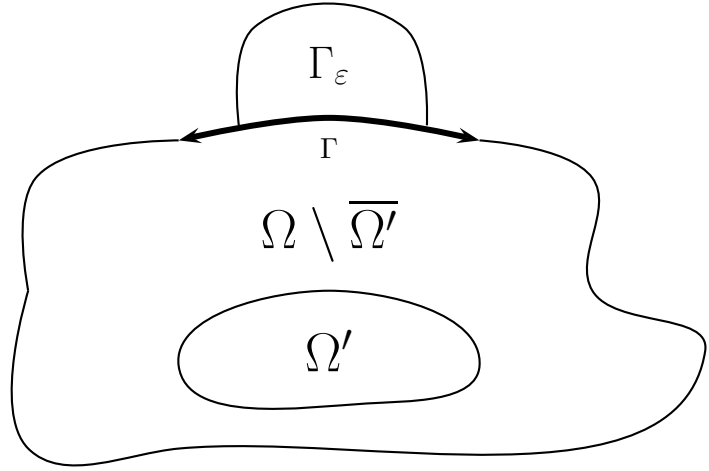

\noindent\textbf{Proof of Lemma \ref{lemma:uniquecont}}. Part i)
of Lemma \ref{lemma:uniquecont} can be proven by trivially
rewriting Theorem 6.5.1's proof in \cite{EK} for the vector case
taking
$W_1=\norm{\widetilde{Q}}_{L^{\infty}(\Omega;\mathcal{M}_{8\times8})}$,
$W_2=0$, $C=\emptyset$. Part ii) follows from part i), remarking
that the boundary conditions on $\Gamma$ guarantee that the
extension of the solution by zero on a neighbourhood
$\Gamma_{\varepsilon}$ in $\R^3$ such that
$\Gamma_{\varepsilon}\cap \Omega=\emptyset$ and
$\overline{\Gamma_{\varepsilon}}\cap \overline{\Omega}$ is an open
subset of $\Gamma$ with respect to the relative topology on
$\partial\Omega$ induced by the Euclidean topology of $\R^3$,
satisfies the same equation and maintains the $H^2$-regularity on
int$((\,\overline{\Omega}\setminus \Omega')\cup
\overline{\Gamma_{\varepsilon}}\,)$. Indeed, the kernel of the
trace operator $(u|_{\partial\mathcal{G}} , (\partial u/\partial
N)|_{\partial\mathcal{G}})$ defined for $u\in H^2(\mathcal{G})$,
is the closure of $C^{\infty}_0(\mathcal{G})$ in
$H^2(\mathcal{G})$ (usually denoted by $H^2_0(\mathcal{G})$), for
any domain $\mathcal{G}$ with $C^{1,1}$ boundary
$\partial\mathcal{G}$ (see \cite{LM} or e.g. \cite[Theorem
1.5.1.5]{G}). For clarity Figure \ref{fig} illustrates the sets
$\Gamma_{\varepsilon}$, $\Omega$, $\Omega'$ in the plane (although
they must be considered in $\R^3$).

\begin{flushright}
$\square$\end{flushright}

\noindent\textbf{Proof of Lemma \ref{lemma:uniquecontL}} Under the
conditions of part i), define $V:=(i\omega\mu)^{-1}\nabla\times U$
and check that $(U,V)$ solves in $\mathcal{G}$ the Maxwell
equations $\nabla\times V +i\omega \gamma U = 0$, $\nabla\times U
-i\omega\mu V = 0$, and $U, V\in H(\mathcal{G};\mbox{curl})$ with
$\nabla\cdot U$, $\nabla\cdot V\in L^2(\mathcal{G})$. By
\cite[Chapter I, Corollary 2.10]{GR}, $U, V\in H^1_{\mbox{\tiny
loc}}(\mathcal{G};\C^3)$. Consider another open ball $B'$ with
$B\cap B'\neq \emptyset$ and $B'\subset \mathcal{G}$. Since the
restrictions of $U,V$ to $B'$ are in $H^1(B';\C^3)$, from the
unique continuation result across $C^2$-surfaces by Eller and
Yamamoto \cite[Corollary 1.2]{EY} for the Maxwell system with
$C^1$ coefficients, we deduce that $U=V=0$ on $B'$. Propagating
this argument we conclude that $U$ and $V$ vanish on any
neighbourhood in $\mathcal{G}$. This proves part i).

Let $\Gamma_{\varepsilon}$ be a nonempty, open, connected subset
of $\R^3$ with Lipschitz boundary such that
$\Gamma_{\varepsilon}\cap \mathcal{G}=\emptyset$,
$\overline{\Gamma_{\varepsilon}}\cap \overline{\mathcal{G}}$ is an
open subset of $\Gamma'$ with respect to the relative topology on
$\partial\mathcal{G}$ induced by the Euclidean topology of $\R^3$.
Figure \ref{fig} with $\Omega=\mathcal{G}$, $\Gamma=\Gamma'$
illustrates the choice of $\Gamma_{\varepsilon}$ in the plane. The
conditions of part ii) guarantee that the extension
$\widetilde{U}$ of $U$ by zero on $\Gamma_{\varepsilon}$ verifies
$\widetilde{U}\in H(\mathcal{G}';\mbox{curl})$, $\nabla\times
(\nabla\times \widetilde{U})\in L^2(\mathcal{G}';\C^3)$ and
$L\,\widetilde{U}=0$ in $\mathcal{G}'$, where
$\mathcal{G}':=\,$int$(\,\overline{\mathcal{G}}\cup\overline{
\Gamma_{\varepsilon}}\,)$. This property follows from the fact
that the $H(\mbox{curl})$-vector functions on a bounded, Lipschitz
domain that can be approximated by smooth compactly supported
functions in $H(\mbox{curl})$-norm are exactly those ones with
zero tangential trace (see e.g. \cite[Theorem 3.33]{Mo} for
details). By part i) of this Lemma \ref{lemma:uniquecontL},
$\widetilde{U}=0$ in $\mathcal{G}'$. In particular, $U=0$ in
$\mathcal{G}$.

\begin{flushright}
$\square$\end{flushright}

Regarding the aforementioned result in \cite{EY}, note that a
counterexample for the stationary Maxwell system with coefficients
in the H\"older class $C^{\alpha}$ for every $\alpha<1$ is
provided in \cite{D} by Demchenko.

\vskip 3mm

\noindent\textbf{Proof of Proposition \ref{prop:density}}.
Following the lines of Lemma 2's proof in \cite{AU}, suppose $v\in
K(\Omega)$ satisfies
$(g|v)_{\Omega'}=\int_{\Omega'}v^{\ast}g\,dx=0$ for any
$g\in\widetilde{K}(\Omega)$. We are going to prove that $v=0$ in
$\Omega$.

Consider the Dirichlet Green's function $G$ in $\Omega$ verifying
for $x\in\Omega$,

$(-\Delta_y I_8) G(x,y) + G(x,y) \widetilde{Q}(y) =
\delta(y-x)I_8$, for $y\in\Omega$, $ G(x,y)=0$, for
$y\in\partial\Omega$, where $\delta$ denotes the Dirac delta
function with pole at the origin and $I_8$ the $8\times 8$
identity matrix. For $g\in\widetilde{K}(\Omega)$ and $x\in\Omega$
we have, by Green's formula,
\begin{gather}
\notag g(x) = \int_{\Omega}(\delta(y-x)I_8) g(y) dy =
\int_{\Omega}(-(\Delta_y I_8)G(x,y) + G(x,y) \widetilde{Q}(y)) g(y) dy\\
\notag = \int_{\Omega} G(x,y) (-\Delta_y I_8 +
\widetilde{Q}(y))g(y)dy - \int_{\partial\Omega}
\bigg(\bigg({\partial \over\partial\nu
(y)}\,I_8\bigg)G(x,y)\bigg)g(y)\bigg)ds(y)\\
\notag = - \int_{\Gamma} \bigg(\bigg({\partial \over\partial\nu
(y)}\,I_8\bigg)G(x,y)\bigg)g(y)\, ds(y).
\end{gather}
In particular note that for $x\in\Omega$,
$U(x)=\int_{\Omega}G(x,y) F(y)\, dy$ provided that $(-\Delta I_8 +
\widetilde{Q})U=F$ in $\Omega$ and $U=0$ on $\partial\Omega$.

By Fubini's theorem,
\begin{gather*}
\int_{\Gamma} \int_{\Omega'} v(x)^{\ast}\bigg(\bigg({\partial
\over\partial\nu (y)}\,I_8\bigg)G(x,y)\bigg) dx \, g(y)\, ds(y)\\
= \int_{\Omega'}v(x)^{\ast}\int_{\Gamma} \bigg(\bigg({\partial
\over\partial\nu (y)}\,I_8\bigg)G(x,y)\bigg) g(y)\, ds(y)\, dx\\
=-\int_{\Omega'}v(x)^{\ast}g(x)\, dx = 0.
\end{gather*}
Thus, for $y\in\Gamma$,
\begin{equation}\label{form:11}
\int_{\Omega'} v(x)^{\ast}\bigg(\bigg({\partial \over\partial\nu
(y)}\,I_8\bigg)G(x,y)\bigg) dx = 0.
\end{equation}

Define the vector field $u$ by
$$
u(y)^{\ast} :=\int_{\Omega'} v(x)^{\ast} G(x,y)\, dx.
$$

\noindent Since $G(x,y)=0$ for $y\in\partial\Omega$ and
$x\in\Omega$, $u(y)=0$ for $y\in\partial\Omega$. By
\eqref{form:11},
\begin{equation}\label{form:12}
\bigg({\partial \over\partial\nu
(y)}\,I_8\bigg)u(y)=0\qquad\mbox{for }y\in\Gamma.
\end{equation}

\noindent Since $(-\Delta I_8 + (\widetilde{Q})^{\ast})u = 0$ in
$\Omega\setminus\overline{\Omega'}$, $u\big|_{\partial\Omega}=0$
and by \eqref{form:12}, it follows that $u=0$ in
$\Omega\setminus\Omega'$ by the unique continuation principle
(Lemma \ref{lemma:uniquecont}). In particular,
$$
u = \bigg({\partial\over\partial\nu} I_8\bigg)u=0\qquad\mbox{on
}\partial\Omega'.
$$

Note that $(-\Delta I_8+(\widetilde{Q})^{\ast})u=v$ in $\Omega'$. Now, we can
write
\begin{gather}
\notag \int_{\Omega'}v(y)^{\ast} v(y)\, dy = \int_{\Omega'}
((-\Delta_y
I_8)u(y)^{\ast} +  u(y)^{\ast} \widetilde{Q}(y))v(y)\, dy\\
\label{form:13}=-\int_{\Omega'} u(y)^{\ast}(\Delta I_8) v(y) dy\\
\notag + \int_{\partial\Omega'} \bigg(u(y)^{\ast}\bigg({\partial
\over
\partial\nu(y)}\,I_8\bigg)v(y)- \bigg(\bigg({\partial
\over\partial\nu (y)}\,I_8\bigg)u(y)^{\ast}\bigg)v(y)\bigg) ds(y)\\
\notag + \int_{\Omega'} u(y)^{\ast} \widetilde{Q}(y) v(y) dy\\
\notag =\int_{\Omega'} u(y)^{\ast} (-\Delta I_8+\widetilde{Q}(y))
v(y) dy=0,
\end{gather}
where identity \eqref{form:13} follows from Green's formula. Hence
$v=0$ in $\Omega'$. Since, $(-\Delta I_8+\widetilde{Q}) v=0$ in
$\Omega$, by unique continuation (Lemma \ref{lemma:uniquecont}),
$v=0$ in $\Omega$.

\begin{flushright}
$\square$\end{flushright}

\noindent\textbf{Proof of Proposition
\ref{prop:density_curlcurl}}. Fix $E\in N(\Omega)$ such that $\int_{\Omega'} \overline{E}\cdot
\widetilde{E}\, dx = 0$ for any $\widetilde{E}\in
\widetilde{N}(\Omega)$. Define $E'$ as the solution to the equation $L\,E' = \chi_{\Omega'}\overline{E}$ in $\Omega$ satisfying $N\times E'=0$ on $\partial\Omega$ in the trace sense, such that $E'\in H(\Omega;\mbox{curl})$ and $\nabla\times\nabla \times E'\in L^2(\Omega;\C^3)$. Using the integration by parts formula \eqref{form:intpartsL}, for any $\widetilde{E}\in
\widetilde{N}(\Omega)$ we have
\begin{gather}
\label{form:41}0=\int_{\Omega'} \overline{E}\cdot \widetilde{E}\, dx = \int_{\Omega} L\, E'\cdot \widetilde{E}\, dx
 = \int_{\Omega} E'\cdot L\,\widetilde{E}\, dx\\ \label{form:42}- \int_{\partial\Omega} {1\over\mu}\, [(N\times E')\cdot (N\times (N\times \nabla\times \widetilde{E})) + (\nabla\times E')\cdot (N\times \widetilde{E})]\, ds\\
 \label{form:43}= -\int_{\Gamma}{1\over\mu}\, (\nabla\times E')\cdot (N\times \widetilde{E})\, ds.
\end{gather}

The trace of $\nabla\times E'$ on $\partial\Omega\in C^{1,1}$ can be decomposed into its tangential and normal components as follows:
\begin{equation}
\label{form:44}\nabla\times E'|_{\partial\Omega} = -N\times (N\times (\nabla\times E'))|_{\partial\Omega} + (N\cdot (\nabla\times E'))\, N|_{\partial\Omega}.
\end{equation}

On using the identity \eqref{form:44} in the integral \eqref{form:43}, the second term in the right hand side of \eqref{form:44} gets cancelled. Therefore from \eqref{form:41}-\eqref{form:43} and \eqref{form:44} we deduce for each $\widetilde{E}\in
\widetilde{N}(\Omega)$,
\begin{equation*}
0= \int_{\Gamma}{1\over\mu}\, N\times (N\times (\nabla\times E'))|_{\partial\Omega}\cdot (N\times \widetilde{E})|_{\partial\Omega}\, ds.
\end{equation*}
So, $N\times (N\times (\nabla\times E'))|_{\partial\Omega}$ vanishes on $\Gamma$. As a result, $N\times (\nabla\times E')|_{\partial\Omega}=0$ on $\Gamma$.

From the condition $N\times E'|_{\partial\Omega}=0$ and the properties $N\times
(\nabla\times E')|_{\partial\Omega}=0$ on $\Gamma$ and $L(E'|_{\Omega\setminus\overline{\Omega '}}) = 0$ in $\Omega\setminus \overline{\Omega'}$, it follows that $E'=0$ in $\Omega\setminus\Omega'$ by the uniqueness result stated in part ii) of Lemma
\ref{lemma:uniquecontL} for $\mathcal{G}=\Omega\setminus \overline{\Omega'}$ and $\Gamma'=\Gamma$. In
particular,
\begin{equation}\label{form:45}
E'|_{\partial\Omega'}=0,\qquad \nabla\times E'|_{\partial\Omega'}
=0.
\end{equation}

Now, we write
\begin{gather*}
\int_{\Omega'} E^{\ast}\, E\, dx = \int_{\Omega'} (L\,
E')^t\, E\, dx
= \int_{\Omega'} (E')^t\, L\, E\, dx\\
-\int_{\partial\Omega'} \mu^{-1} ((E')^t\, (N\times
(\nabla\times E)) + (\nabla\times E')^t\, (N\times E)
)\, ds = 0,
\end{gather*}
where the last identity follows from formula \eqref{form:intpartsL}.
Hence, $E=0$ in $\Omega'$. Since $L\, E=0$ in $\Omega$, we deduce by
the unique continuation principle stated in part i) of Lemma
\ref{lemma:uniquecontL} with $\mathcal{G}=\Omega$, that $E=0$ in $\Omega$.

\begin{flushright}
$\square$\end{flushright}

\section{Proof of uniqueness}

Here, the outline of \cite{A}, Section 3.2 in \cite{Ca10}, Section
3.4 in \cite{Ca11} is adapted to prove Theorem \ref{theo:main}.

Let $\omega>0$ be the time-harmonic frequency. Assume
$\mu_j$,$\gamma_j$ is an admissible pair of coefficients for each
$j=1,2$, according to Definition \ref{def:admissible}, such that
$\supp(\mu_1-\mu_2)$, $\supp (\gamma_1-\gamma_2)\subset\Omega$.
Let $\Omega'$ be an open subset of $\R^3$ with $C^2$ boundary such
that $\Omega'\subset\subset\Omega$,
$\Omega\setminus\overline{\Omega'}$ is connected and
\begin{equation}\label{form:equalityneighbour}
\mu_1=\mu_2\qquad\mbox{and}\qquad\gamma_1=\gamma_2\qquad\mbox{in
}\overline{\Omega}\setminus\Omega'.
\end{equation}
Suppose $C_{\Gamma}^1 = C_{\Gamma}^2$. The extended coefficients
to $\R^3$ according to the extensions described in Section
\ref{section:solutions} will be written likewise,
$\mu_j$,$\gamma_j$. Remember that $\mu_j=\mu_0$,
$\gamma_j=\varepsilon_0$ outside the ball $B(\mbox{O},\rho)$,
where $\overline{\Omega}\subset B(\mbox{O},\rho)$, and $\mu_0$,
$\varepsilon_0$ are constants.

Let $j\in\{1,2\}$. Define
\begin{gather*}
\alpha_j :=\log\gamma_j,\qquad \beta_j:=\log\mu_j,\qquad
\kappa_j:=\omega\mu_j^{1/2}\gamma_j^{1/2},\\
f := \chi_{\Omega}\cdot \bigg({1\over 2}\Delta (\alpha_1-\alpha_2)
+ {1\over 4} \,(\nabla\alpha_1\cdot\nabla\alpha_1 -
\nabla\alpha_2\cdot\nabla\alpha_2) +
(\kappa_2^2-\kappa_1^2)\bigg),\\
g := \chi_{\Omega}\cdot \bigg({1\over 2}\Delta (\beta_1-\beta_2) +
{1\over 4} \,(\nabla\beta_1\cdot\nabla\beta_1 -
\nabla\beta_2\cdot\nabla\beta_2) + (\kappa_2^2-\kappa_1^2)\bigg),
\end{gather*}
where $\chi_{\Omega}$ denotes the characteristic function of
$\Omega$.

Fix $\xi\in\R^3\setminus 0$. Let $\S^2$ denote the unit sphere in
$\R^3$. Assume $\tau\geq 1$ and take $\eta_1,\eta_2\in\S^2$ with $
\eta_2\cdot\eta_1 = \eta_1\cdot\xi = \eta_2\cdot\xi = 0 $, and
\begin{gather*}
\zeta_1 = -{1\over 2}\,\xi + i\bigg( \tau^2 + {|\xi|^2\over 4}
\bigg)^{1/2}\eta_1 + (\tau^2 +
\omega^2\varepsilon_0\mu_0)^{1/2}\,\eta_2,\\
\zeta_2 = {1\over 2}\,\xi - i\bigg( \tau^2 + {|\xi|^2\over 4}
\bigg)^{1/2}\eta_1 + (\tau^2 +
\omega^2\varepsilon_0\mu_0)^{1/2}\,\eta_2.
\end{gather*}

Note that $\zeta_j\in\C^3$ satisfies $\zeta_j\cdot\zeta_j =
\omega^2\varepsilon_0\mu_0$, and
$$
|\zeta_j| = (|\mbox{Re}(\zeta_j)|^2 +
|\mbox{Im}(\zeta_j)|^2)^{1/2} = (|\xi|^2/2 + 2\tau^2 +
\omega^2\varepsilon_0\mu_0)^{1/2}.
$$
Further, $ \zeta_1-\overline{\zeta}_2 = -\xi$, and as
$\tau\rightarrow\infty$,
$$
{\zeta_1\over|\zeta_1|} = {1\over\sqrt{2}}\, (i\,\eta_1+\eta_2) +
\mathcal{O}(\tau^{-1}),\qquad {\zeta_2\over|\zeta_2|} =
{1\over\sqrt{2}}\, (-i\,\eta_1+\eta_2) + \mathcal{O}(\tau^{-1}),
$$
where the implicit constants depend on $|\xi|$ (and
$\omega$,$\varepsilon_0$,$\mu_0$).

Consider the special solutions $$Z_1(x,\zeta_1)=e^{i\zeta_1\cdot
x}(L_1(\zeta_1)+R_1(x,\zeta_1)),\qquad Y_2(x,\zeta_2) =
e^{i\zeta_2\cdot x}(M_2(\zeta_2)+S_2(x,\zeta_2))$$from Proposition
\ref{prop:special1} and Proposition \ref{prop:special2} applied to
the case $\mu=\mu_1$, $\gamma=\gamma_1$, $\zeta=\zeta_1$ and
$\mu=\mu_2$, $\gamma=\gamma_2$, $\zeta=\zeta_2$, respectively, so
that $Z_1$,$Y_2$ solve $(-\Delta I_8 +Q_1)Z_1= 0$,
$(P+W_2^{\ast})Y_2=0$ in $\R^3$. Choosing such solutions with
$B_j=0$ and $A_j$ such that
$$
\bigg(i{\eta_1\over \sqrt{2}} + {\eta_2\over \sqrt{2}}\bigg)\cdot
A_1 = \bigg(i{\eta_1\over \sqrt{2}} + {\eta_2\over
\sqrt{2}}\bigg)\cdot \overline{A_2} = 1,
$$
one obtains
\begin{equation}\label{form:30}
((Q_1-Q_2)Z_1|Y_2)_{\Omega} = \hat{f}(\xi) +
\mathcal{O}(\tau^{-1}),
\end{equation}
as $\tau\rightarrow\infty$. Analogously, choosing $Z_1$, $Y_2$
with $A_j=0$ and $B_j$ such that
$$
\bigg(i{\eta_1\over \sqrt{2}} + {\eta_2\over \sqrt{2}}\bigg)\cdot
B_1 = \bigg(i{\eta_1\over \sqrt{2}} + {\eta_2\over
\sqrt{2}}\bigg)\cdot \overline{B_2} = 1,
$$
one can prove
\begin{equation}\label{form:31}
((Q_1-Q_2)Z_1|Y_2)_{\Omega} = \hat{g}(\xi) +
\mathcal{O}(\tau^{-1}),
\end{equation}
as $\tau\rightarrow\infty$. In \eqref{form:30},\eqref{form:31} the
implicit constant depends on
$M,\xi,|\Omega|,\rho,\omega,\varepsilon_0,\mu_0$.

Fix $\epsilon>0$. For each choice of $Z_1$, define
$Y_1:=(P-W_1^t)Z_1$. Hence, $(P+W_1)Y_1=0$ in $\R^3$ and by
Proposition \ref{prop:special1}, $Y_1$ reads $ Y_1 =
(0\,\,\,\,\,\mu_1^{1/2}
H_1^t\,|\,0\,\,\,\,\,\gamma_1^{1/2}E_1^t)^t, $ with $E_1,H_1$
solutions of
\begin{gather}
\notag \nabla\times H_1 +i\omega \gamma_1 E_1 = 0,\\
\label{form:40}\nabla\times E_1 -i\omega\mu_1 H_1 = 0,
\end{gather}
in $\R^3$. In particular, $L_1\, E_1 = 0$ in $\Omega$, where $L_1$
denotes the second order operator $L$ defined in
\eqref{form:secondorder} for $\mu_1$,$\gamma_1$. By Proposition
\ref{prop:special1}, $Z_1|_{\Omega}\in H^2(\Omega;\C^8)$. Thus, by
the Lipschitz regularity and the a priori bounds from below for
$\mu_1$, $\gamma_1$, and from equation \eqref{form:40} we deduce
that $E_1|_{\Omega}\in H(\Omega;\mbox{curl})$ and $\nabla\times
(\nabla\times E_1)|_{\Omega}\in L^2(\Omega;\C^3)$.

By Proposition \ref{prop:density_curlcurl}, there exists
$\widetilde{E}_1\in H(\Omega;\mbox{curl})$ such that $\nabla\times
(\nabla\times \widetilde{E}_1)\in L^2(\Omega;\C^3)$, $L_1\,
\widetilde{E}_1=0$ in $\Omega$, $N\times
\widetilde{E}_1|_{\partial\Omega}=0$ on $\Gamma_c$, and $
\norm{E_1-\widetilde{E}_1}_{L^2(\Omega';\C^3)}<\epsilon$.

Due to the a priori condition $\mu_1\geq M^{-1}$, $\mu_1$ does not
vanish. Define $\widetilde{H}_1:=(1/i\omega\mu_1)\nabla\times
\widetilde{E}_1$. Therefore $\widetilde{E}_1$,$\widetilde{H}_1\in
H(\Omega;\mbox{curl})$ solve
\begin{gather*}
\nabla\times \widetilde{H}_1 +i\omega \gamma_1 \widetilde{E}_1 = 0,\\
\nabla\times \widetilde{E}_1 -i\omega\mu_1 \widetilde{H}_1 = 0,
\end{gather*}
in $\Omega$. Define $ \widetilde{Y}_1 := (0\,\,\,\,\,\mu_1^{1/2}
\widetilde{H}_1^t\,|\,0\,\,\,\,\,\gamma_1^{1/2}\widetilde{E}_1^t)^t$.

\subsection{Invertibility of the Dirac-type operator $P$}\label{subsection}

$\,$

Along this subsection the letters $M$, $N$, $E$ refer to mathematical entities which are different from their meanings in the rest of the paper.

Jochen Br\"uning and Matthias Lesch in \cite{BL} generalize the
analysis of Dirac-type operators considered in the well-known
paper by Atiyah, Patodi and Singer \cite{APS}. Concerning the
general Dirac-type operators studied there on compact manifolds
with boundary, in \cite[Section 1.B]{BL} an operator $D$ is
introduced acting on sections of a hermitian vector bundle $E$
over an open subset $M$ of a compact oriented Riemannian manifold
$\widetilde{M}$ such that its boundary $N=\partial M$ is a compact
hypersurface in $\widetilde{M}$. The authors call $\widetilde{E}$
the vector bundle over $\widetilde{M}$, and $E_N:=
\widetilde{E}\upharpoonright N$. The differential operator $D$ is
said to be of Dirac type if it is first order, symmetric and
elliptic in $L^2(E)$ with domain $C^{\infty}_0(E)$ verifying that
$D^2$ has scalar principal symbol given by the metric tensor.

Taking $M=\Omega'$ and $E=M\times \C^8$ the trivial bundle over $M$, each fiber equipped with the standard hermitian inner product of $\C^8$, the operator $P$ on $\Omega'$ defined in \eqref{form:operatorP} falls into the category of these Dirac type operators, since $P^2 =-\Delta I_8$, $\langle PU, V \rangle_{\Omega'} = \langle U, P V \rangle_{\Omega'}$ for any $U$, $V\in C^{\infty}_0(\Omega'; \C^8)$ and the characteristic form of $P$, namely $Q(\lambda)=\mbox{det}(\Lambda(\lambda))=-i|\lambda|^8$, does not vanish for any $\lambda = (\lambda_1,\lambda_2,\lambda_3)\in\R^3\setminus 0$. Here, $\Lambda(\lambda)$ denotes the symbol of $P$ given by the matrix form
$$
\Lambda(\lambda)={1\over i}  \left(%
\begin{array}{c|c}
  0 & \mathcal{A}(\lambda) \\
  \hline
  \mathcal{B}(\lambda) & 0  \\
\end{array}%
\right)\begin{array}{c}
\phantom{.}\\
,\\
\end{array}
$$
with
$$
 \mathcal{A}(\lambda)= \left( \begin{array}{cccc}
  0 & \lambda_1 & \lambda_2 & \lambda_3 \\
  \lambda_1 & 0 & \lambda_3 & -\lambda_2\\
  \lambda_2 & -\lambda_3 & 0 & \lambda_1\\
  \lambda_3 & \lambda_2 & -\lambda_1 & 0\\
\end{array}%
\right)\begin{array}{c}
\phantom{.}\\
\phantom{.}\\
\phantom{.}\\
,\\
\end{array}\qquad
\mathcal{B}(\lambda)= \left( \begin{array}{cccc}
  0 & \lambda_1 & \lambda_2 & \lambda_3 \\
  \lambda_1 & 0 & -\lambda_3 & \lambda_2\\
  \lambda_2 & \lambda_3 & 0 & -\lambda_1\\
  \lambda_3 & -\lambda_2 & \lambda_1 & 0\\
\end{array}%
\right)\begin{array}{c}
\phantom{.}\\
\phantom{.}\\
\phantom{.}\\
.\\
\end{array}
$$

In \cite{BL} it is proved that $D$ admits self-adjoint extensions
by imposing non-local boundary conditions given by an orthogonal
projection $\pi$ in $L^2(E_N)$, which is a classical
pseudodifferential operator on $E_N$ satisfying a certain symmetry
property (condition (1.13) in \cite{BL}) related to the structure
of the operator (see \cite[Lemma 1.1]{BL} for details). For such
$\pi$ and by \cite[Theorem 1.5]{BL} and the interpretation by Y.
Kurylev and M. Lassas \cite[Theorem 2.1]{KL}, it turns out that
$P\upharpoonright \mathcal{D}$ is self-adjoint with empty
essential spectrum and finite-dimensional eigenspaces, where
$\mathcal{D}:= \{ U\in H^1(\,\overline{\Omega'}\, ; \C^8)\, :\,
\pi(U|_{\partial\Omega'})=0 \}$.

The domain $\mathcal{D}$ in $H^1(\Omega';\C^8)$ with the graph
norm associated with $P$ is continuously embedded into
$H^1(\Omega' ; \C^8)$. The space $H^1(\Omega' ; \C^8)$ is
compactly embedded into $L^2(\Omega' ; \C^8)$. Since $W\in
L^{\infty}(\Omega' ; \mathcal{M}_{8\times8})$ for admissible
$\mu$, $\gamma$, the operator of multiplication by $W^t$, which we
write $\mathcal{M}_{W^t}$, is bounded and linear in $L^2(\Omega' ;
\C^8)$. Therefore, $\mathcal{M}_{W^t}$ is $(P\upharpoonright
\mathcal{D})$-compact.

Thus, $P-W^t$ has also empty essential spectrum and
finite-dimensional ei\-gen\-spa\-ces. If $0$ is in the spectrum of
$P-W^t$, then $0$ must be an eigenvalue with finitely many
linearly independent eigen- and associated functions. We can make
$0$ no longer be an eigenvalue by choosing a new set of boundary
conditions which are not satisfied by any of the finitely many
linearly independent eigen- and associated functions in the root
spaces associated with $0$. Let us keep denoting the resultant
boundary operator by $\pi$ so that the condition $\pi
(Z|_{\partial\Omega'})=0$ guarantees the existence of a constant
$C_{\mbox{\tiny stblty}}$ independent of $Z$ such that
\begin{equation}\label{form:stabltydirac}
\norm{Z}_{L^2(\Omega';\C^8)}\leq C_{\mbox{\tiny stblty}} \norm{(P-W^t)
Z}_{L^2(\Omega';\C^8)},
\end{equation}
provided that $Z, (P-W^t) Z \in L^2(\Omega';\C^8)$.

\vskip 5mm

The argument presented in Subsection \ref{subsection}, together
with a trick based on an auxiliary system which improves the
regularity of $(P-W^t)Z$ when Maxwell equations are satisfied,
leads to the following
\begin{lemma}\label{lemma:dirac}
For admissible coefficients $\mu$,$\gamma$, assume that $(P-W^t)Z
= Y$ in $\Omega'$, where $Y$ reads $ Y = (0\,\,\,\,\,\mu^{1/2}
H^t\,|\,0\,\,\,\,\,\gamma^{1/2}E^t)^t, $ with $E,H$ verifying
\eqref{form:maxwelleq} in $\Omega'$. Additionally, suppose $\pi(Z|_{\partial\Omega'})=0$ for the boundary operator $\pi$ introduced in Subsection \ref{subsection}. Then there exists a constant
$C$ only depending on $C_{\mbox{\tiny stblty}}$, $M$, $\omega$,
such that $\norm{Z}_{L^2(\Omega';\C^8)}\leq C\,
\norm{E}_{L^2(\Omega';\C^3)}$.
\end{lemma}

\noindent\textbf{Proof of Lemma \ref{lemma:dirac}}. Writing $Z =
(Z_1\,\,\,\,\,Z_H^t\,|\,Z_2\,\,\,\,\,Z_E^t)^t$, where $Z_j$
($j=1,2$) are scalar fields and $Z_H$, $Z_E$ vector fields, and
defining $Z_{\mbox{{\tiny
aux}}}:=(Z_1\,\,\,\,\,Z_H^t\,|\,Z_2\,\,\,\,\,(Z_{E}')^t)^t$ with
$Z_{E}':= Z_E + \mu^{-1/2}\omega^{-1}E$, the dependence on the
electric field $E$ of the vector function $(P-W^t)Z_{\mbox{\tiny
aux}} $ is zeroth-order. Indeed, it is straightforward to check
that
$$
(P-W^t)Z_{\mbox{\tiny aux}} = \left(
\begin{array}{c}
(-i/\omega) (-2\nabla\mu^{-1/2}+\mu^{-1/2}\nabla\alpha)\cdot E\\
0\\
\hline
0\\
(\gamma^{1/2}-\kappa\mu^{-1/2}\omega^{-1})E\\
\end{array}
\right)\begin{array}{c}
\phantom{.}\\
\phantom{.}\\
\phantom{.}\\
.\\
\end{array}
$$

By \eqref{form:stabltydirac},
\begin{align*}
\norm{Z}_{L^2(\Omega';\C^8)}&\leq \norm{Z_{\mbox{\tiny
aux}}}_{L^2(\Omega';\C^8)} +
M^{1/2}\omega^{-1}\norm{E}_{L^2(\Omega';\C^3)}\\
&\leq C_{\mbox{\tiny stblty}}\norm{(P-W^t)Z_{\mbox{\tiny
aux}}}_{L^2(\Omega';\C^8)} +
M^{1/2}\omega^{-1}\norm{E}_{L^2(\Omega';\C^3)}\\
&\leq (C_{\mbox{\tiny stblty}} C(M,\omega) +
M^{1/2}\omega^{-1})\norm{E}_{L^2(\Omega';\C^3)}.
\end{align*}

\begin{flushright}
$\square$\end{flushright}

For $\widetilde{Y}_1$ defined above, let $\widetilde{Z}_1$ be a solution to the system $(P-W_1^t)
\widetilde{Z}_1 = \widetilde{Y}_1$ in $\Omega'$ such that $\pi((Z_1-\widetilde{Z}_1)|_{\partial\Omega'})=0$. By Lemma \ref{lemma:dirac},
\begin{gather}\label{form:28}
\norm{Z_1-\widetilde{Z}_1}_{L^2(\Omega';\C^8)}\leq C\,
\norm{E_1-\widetilde{E}_1}_{L^2(\Omega';\C^3)}\leq \epsilon C.
\end{gather}

For each choice of $Y_2$ (with $Y_2|_{\Omega}\in H^1(\Omega;\C^8)$), since $(-\Delta I_8 +
\widetilde{Q}_2\, )Y_2=0$ in $\Omega$ by Lemma \ref{lemma:schro},
where $\widetilde{Q}_2$ denotes the zeroth order matrix operator
$\widetilde{Q}$ defined in \eqref{form:matpotentilde} for
$\mu_2$,$\gamma_2$, by elliptic regularity and Proposition \ref{prop:density} there exists
$\widetilde{Y}_2\in H^2(\Omega;\C^8)$ verifying $(-\Delta I_8 +
\widetilde{Q}_2\, )\widetilde{Y}_2=0$ in $\Omega$ with
$\widetilde{Y}_2|_{\partial\Omega}=0$ on $\Gamma_c$, and
\begin{equation}\label{form:29}
\norm{Y_2-\widetilde{Y}_2}_{L^2(\Omega';\C^8)}<\epsilon.
\end{equation}

Condition \eqref{form:equalityneighbour} implies that
$\mu_1=\mu_2$, $\gamma_1=\gamma_2$,
$\partial_{x_l}\mu_1=\partial_{x_l}\mu_2$ and
$\partial_{x_l}\gamma_1=\partial_{x_l}\gamma_2$ on
$\overline{\Gamma}$ (for $l=1,2,3$). By Proposition
\ref{prop:intboundary},
\begin{gather}\label{form:27}
((Q_1-Q_2)\widetilde{Z}_1|\widetilde{Y}_2)_{\Omega}=0.
\end{gather}

Applying \eqref{form:28} and \eqref{form:29} write
\begin{gather*}
|((Q_1-Q_2)Z_1|Y_2)_{\Omega'}-((Q_1-Q_2)\widetilde{Z}_1|\widetilde{Y}_2)_{\Omega'}
|\\
=\bigg|\int_{\Omega'} (Q_1-Q_2)Z_1\cdot
\bigg(\,\overline{Y_2-\widetilde{Y}_2}\,\bigg)\, dx  -
\int_{\Omega'}
(Q_1-Q_2)(\widetilde{Z}_1-Z_1)\cdot\overline{\widetilde{Y}_2}\,
dx\bigg|\\
\leq \norm{Q_1-Q_2}_{L^{\infty}(\Omega';\mathcal{M}_{8\times 8})}\,
\bigg(\norm{Z_1}_{L^2(\Omega';\C^8)}\norm{Y_2-\widetilde{Y}_2}_{L^2(\Omega';\C^8)}
\\
+\norm{\widetilde{Y}_2}_{L^2(\Omega';\C^8)}\norm{Z_1-\widetilde{Z}_1}_{L^2(\Omega';\C^8)}\bigg)\\
\leq C(C_{\mbox{\tiny stblty}},M,\omega)\epsilon \,
\bigg(\norm{Z_1}_{L^2(\Omega;\C^8)} +
\norm{\widetilde{Y}_2}_{L^2(\Omega';\C^8)}\bigg)\\
\leq \epsilon
C(M,\Omega,\rho,\xi,\omega,\varepsilon_0,\mu_0,C_{\mbox{\tiny
stblty}})\, e^{c(\tau + |\xi|)},
\end{gather*}
where $c=c(\Omega)$ and last inequality follows from the fact
$$
\norm{\widetilde{Y}_2}_{L^2(\Omega';\C^8)}\leq
\norm{Y_2-\widetilde{Y}_2}_{L^2(\Omega';\C^8)} +
\norm{Y_2}_{L^2(\Omega';\C^8)}\leq 1+ \norm{Y_2}_{L^2(\Omega;\C^8)},
$$
the exponential behavior of $Z_1$, $Y_2$ and estimates
\eqref{form:resR}, \eqref{form:resS}. Therefore, denoting
$$
\epsilon'(\epsilon):= \epsilon
C(\Omega,\rho,\xi,\omega,\varepsilon_0,\mu_0,C_{\mbox{\tiny
stblty}}),
$$
since $Q_1=Q_2$ in $\Omega\setminus\Omega'$ and by
\eqref{form:27}, we have
\begin{gather}
\label{form:32bis}
|((Q_1-Q_2)Z_1|Y_2)_{\Omega}|=|((Q_1-Q_2)Z_1|Y_2)_{\Omega'}|\\
\label{form:32} \leq
|((Q_1-Q_2)\widetilde{Z}_1|\widetilde{Y}_2)_{\Omega'}| +
\epsilon'(\epsilon)\, e^{c(\tau + |\xi|)} = \epsilon'(\epsilon)\,
e^{c(\tau + |\xi|)}.
\end{gather}
Thus, for fixed $\tau$ and $\xi$, letting $\epsilon\rightarrow 0$
in \eqref{form:32bis}-\eqref{form:32}, we get
\begin{equation}\label{form:33}
((Q_1-Q_2)Z_1|Y_2)_{\Omega}=0,
\end{equation}
for both choices of $Z_1$, $Y_2$. By \eqref{form:30},
\eqref{form:31} and \eqref{form:33}, we have for large enough
$\tau$,
$$
|\hat{f}(\xi)| + |\hat{g}(\xi)|\leq {C\over \tau},
$$
where $C=C(M,\xi,|\Omega|,\rho,\omega,\mu_0,\varepsilon_0)$. For
any fixed $\xi\in\R^3$, by letting $\tau\rightarrow\infty$ deduce
that $\hat{f}(\xi)=\hat{g}(\xi)=0$. Hence, $f=g=0$.

Using a Carleman estimate, Pedro Caro in \cite{Ca11} proves the
following inequality
\begin{gather}
\label{form:carleman1}e^{d_1/h}\sum_{j=1,2}
(h\norm{\phi_j}_{L^2(\Omega)}^2 +
h^3\norm{\nabla\phi_j}_{L^2(\Omega)}^2)\leq C\, e^{d_2/h}\\
\label{form:carleman2}\times\bigg( h^4
\big(\norm{f}_{L^2(\Omega)}^2 + \norm{g}_{L^2(\Omega)}^2\big) +
\sum_{j=1,2} \big(h\norm{\phi_j}_{L^2(\partial\Omega)}^2 + h^3
\norm{\nabla\phi_j}_{L^2(\partial\Omega)}^2\big) \bigg),
\end{gather}
where $\phi_1:=\gamma_1^{1/2}-\gamma_2^{1/2}$, $\phi_2 :=
\mu_1^{1/2}-\mu_2^{1/2}$, $C=C(\Omega,M)$, $0<h<C^{-1/3}\leq 1$,
and
$$d_1:=\inf\{|x-x_0|^2\, :\, x\in\Omega\},\qquad d_2:=\sup\{|x-x_0|^2\, :\,
x\in\Omega\},$$for certain point $x_0\notin\overline{\Omega}$.
Under Theorem \ref{theo:main}'s conditions the summation term of
norms on $\partial\Omega$ in \eqref{form:carleman2} vanishes. From
this fact together with $f=g=0$ in $\Omega$, we conclude
$\mu_1=\mu_2$ and $\gamma_1=\gamma_2$ in $\Omega$.

%%%%%%%%%%%%%%%%%%%%%%%%%%%%%%%%%%%%%%%%%%%%%%%%%%%%%%%%%%%%%%%%%%%%%%%%%%%%%%%%%%

%\bibliographystyle{siam}
%\bibliography{Inverse_problems_references}

\end{document}